\documentclass[onefignum,onetabnum]{siamart190516}
\usepackage{amssymb,bbm,amsmath,color,psfrag,subfigure}
\usepackage{graphicx}
\usepackage{stackrel,hyperref}
\usepackage{algorithm,algorithmic} 
\usepackage{wrapfig}
\usepackage{mathtools}
\usepackage{tikz,pgfplots}

\def\qed{\hfill \vrule height 7pt width 7pt depth 0pt\medskip}
\def\beq{\begin{equation}}
\def\eeq{\end{equation}}
\def\proof{\noindent{\bf Proof}\ \ }

\newcommand{\ds}{\displaystyle}

\newcommand{\ba}{\begin{array}}
\newcommand{\ea}{\end{array}}

\newcommand{\be}{\begin{equation}}
\newcommand{\ee}{\end{equation}}

\newcommand{\1}{\mathbbm{1}}
\newcommand{\E}{\mc{E}}

\newcommand{\R}{\mathbb{R}}

\newcommand{\de}{\mathrm{d}}

\DeclareMathOperator{\diag}{diag}

\newcommand{\mc}{\mathcal}

\newcommand{\norm}[1]{\left\lVert#1\right\rVert}

\definecolor{mycolor1}{RGB}{230,97,1}
\definecolor{mycolor2}{RGB}{253,184,99}
\definecolor{mycolor3}{RGB}{178,171,210}
\definecolor{mycolor4}{RGB}{94,60,153}

\newsiamremark{remark}{Remark}
\newsiamremark{example}{Example}
\crefname{hypothesis}{Hypothesis}{Hypotheses}
\newsiamthm{claim}{Claim}

\headers{Well-posedness of Feedback-controlled Dynamical Flow Networks}{G. Como and G. Nilsson}

\title{On the well-posedness of dynamical flow networks with feedback-controlled outflows \thanks{Submitted to the editors DATE.
\funding{This research was carried on within the framework of the MIUR-funded {\it Progetto di Eccellenza} of the {\it Dipartimento di Scienze Matematiche G.L.~Lagrange}, CUP: E11G18000350001, and was partly supported by the {\it Compagnia di San Paolo} and the Swedish Research Council.}}}

\author{Giacomo Como\thanks{Dipartimento di Scienze Matematiche, Politecnico di Torino, Corso Duca degli Abruzzi 24, 10129, Torino, Italy 
  and Department of Automatic Control, Lund University, Sweden  (\email{giacomo.como@polito.it}, \url{http://staff.polito.it/giacomo.como/}).}
\and Gustav Nilsson\thanks{School of Electrical and Computer Engineering, Georgia Institute of Technology, 777 Atlantic Drive NW, Atlanta, GA, USA 
  (\email{gustav.nilsson@gatech.edu}, \url{http://gustavnilsson.name}).}
}

\ifpdf
\hypersetup{
  pdftitle={Well-posedness of dynamical flow networks with feedback-controlled outflows},
  pdfauthor={G. Como and G. Nilsson}
}
\fi

\begin{document}

\maketitle
\begin{abstract}
We study the well-posedness of a class of dynamical flow network systems describing the dynamical mass balance among a finite number of cells exchanging flow of a commodity between themselves and with the external environment. Systems in the considered class are described as differential inclusions whereby the routing matrix is constant and the outflow from each cell in the network is limited by a control that is a Lipschitz continuous function of the state of the network. In many applications, such as queueing systems and traffic signal control, it is common that an empty queue can be allowed to have more outflow than the mass in the queue. While models for this scenario have previously been presented for open-loop outflow controls, this result ensures the existence and uniqueness of solutions for the network flow dynamics in the case Lipschitz continuous feedback controllers. 

\end{abstract}



\section{Introduction}\label{sec:introduction}
The study of dynamical flow network systems has recently gained a good deal of attention, motivated by applications to transportation systems~\cite{coogan2015compartmental, Como.Lovisari.ea:TCNS15} and queueing networks~\cite{kelly2004fluid,massoulie2007structural,walton2014concave}. Such systems describe the dynamical flow of mass among a finite set of interconnected cells, and have sometimes been referred to as compartmental systems in the control literature \cite{Jacquez.Simon:1993,Walter.Contreras:1999}. 

In order to capture congestion effects, such dynamical flow network systems are typically nonlinear \cite{Como:2017}. In particular, most of them prescribe that the outflow from a cell in the network is limited by a nonlinear function of the traffic volumes. For example, in a traffic network, the outflow from some lanes are limited by a traffic signal, and in a queuing network, the service rate is limited by a server's capacity, which may be state-dependent, and possibly an admission controller. Quite common for many of those applications is that there is a feedback controller actively limiting the outflow from each link. However, it is not always the case that there is enough mass on the link to achieve the outflow limit imposed by the controller. In this cases, first-order ODE-based models fall short of describing the network flow dynamics while guaranteeing physically meaningful properties such as  mass conservation and non-negativity of traffic volumes.

In this paper, we study dynamical models for flow network systems, where the actual outflow from the links may be less than what the controller allows for, in the case when there is not enough mass present on the links. The dynamics is then described by a differential inclusion, and under the assumption that the outflow from the links is the maximum allowed when there is mass present on the links, we show that there exists a unique solution to the dynamics. The dynamical model we are using is a point-queue model. This model is sometimes referred to as vertical queues, to emphasize that a possible spatial distribution of the particles queuing up is not considered in the model. The contribution of this paper is that we show the existence and uniqueness of a solution to a dynamical model for flow networks can be guaranteed when the controller is feedback-based and Lipschitz continuous.

A similar point-queue model, but where the outflow controller does not have feedback, i.e., it is an open-loop controller, has been studied for traffic signal control in~\cite{muralidharan2015analysis} and~\cite{hosseini2017}. In those papers, the existence and uniqueness of a solution to the dynamics have been shown when control action is binary, i.e., the traffic signal at a given time point is either green or red. However, in many feedback based solutions for flow network control, one is instead considering an averaged control signal that depends continuously on the state, such as in~\cite{grandinetti2018distributed, bianchin2019gramian, nilsson2017generalized}. In~\cite{canudas2019average} it has been shown that under certain assumptions, the averaged control signal dynamics stays close to the binary control signal dynamics.

\subsection{Notation}
We let $\R_{(+)}$ denote the (non-negative) reals. For a finite set $\mc A$, we let $\R^{\mc A}$ denote the set of vectors with real entries indexed by the elements of $\mc A$. For a vector $a \in \R^n$, we let $\diag(a) \in \R^{n \times n}$ be a matrix with the entries of $a$ on the diagonal and all off-diagonal entires equal to zero. With $\1$ we denote a vector whose all entries equals one. Inequalities between vectors are meant to hold entry-wise, i.e., e.g., $a\le b$ for $a,b\in\R^{\mc A}$ means that $a_i\le b_i$ for every $i\in\mc A$. The positive part and the negative part of a vector $x\in\R^{\mc A}$ are denoted by $[x]_+ = \max(x, 0)\in\R^{\mc A}$ and $[x]_- = \max(- x, 0)\in\R^{\mc A}$, respectively, where $\max$ and $\min$ are applied entry-wise. Analogously, the absolute value of a vector $x\in\R^{\mc A}$ is the vector $|x|=[x]_++[x]_-\in\R_+^{\mc A}$ whose entries are equal to the absolute values of the entries of $x$. We let $\norm{\cdot}$ be the standard $2$-norm, unless otherwise specified. 
Finally, a directed multigraph is a $4$-tuple $\mc G = (\mc V, \mc E,\sigma,\tau)$ where $\mc V$ is a finite set of nodes, $\mc E$ is a finite set of links, and $\sigma,\tau:\mc E\to\mc V$ are the maps assigning to each link $i\in\mc E$ its tail node $\sigma(i)$ and head node $\tau(i)$, respectively, such that $\sigma(i)\ne\tau(i)$ for every $i\in\mc E$. 

\section{Model}\label{sec:model} 
In this section we present the dynamical flow network system. 

We model the network topology as a directed multigraph $\mc G = (\mc V, \mc E,\sigma,\tau)$. Every link $i\in\mc E$ is to be interpreted as a cell containing a traffic volume $x_i=x_i(t)\ge0$, for $t\ge0$. The state of the system is described by the vector $x=x(t)\in\mc X$ where $\mc X \in \R^{\mc E}_+$ whose entries coincide with the traffic volumes in the different links and evolves in continuous time as the cells exchange flow with adjacent cells and possibly with the external environment.  

In particular, let each cell $i\in\mc E$ possibly receive an exogenous inflow $\lambda_i=\lambda_i(t)\ge0$ directly from the external environment. Moreover, let $z_i=z_i(t)\ge0$ be the total outflow from cell $i$ directed towards immediately downstream cells and possibly to the external environment. Specifically, we shall assume that a constant fraction $R_{ij}\ge0$ of the outflow $z_i$ from cell $i$ flows towards a cell $j\ne i$ such that $\tau(i)=\sigma(j)$, while the remaining part $(1-\sum_jR_{ij})z_i$ leaves the network directly. Conservation of mass then gives 
\be\label{eq:dyn1}\dot x_i=\lambda_i+\sum_{\substack{j\in\mc E:\\\tau(j)=\sigma(i)}}R_{ji}z_j-z_i\,,\qquad i\in\mc E\,.\ee

In order to introduce a more compact notation, let us stack the exogenous inflows in a vector $\lambda=(\lambda_i)_{i\in\mc E}\in\R_+^{\mc E}$ and the cells' outflows in a vector $z=(z_i)_{i\in\mc E} \in\R_+^{\mc E}$. Moreover, let us introduce the routing matrix $R \in \R_+^{\E \times \E}$ whose entries $R_{ij} \geq 0$ coincide with fraction of outflow from cell $i$ that flows directly to cell $j$. Observe that the network topology constraints imply that $R_{ij}=0$ whenever $\tau(i)\ne\sigma(j)$. Moreover, conservation of mass implies that $\sum_j R_{ij}\le 1$ for every cell $i$, i.e., the routing matrix $R$ has be row sub-stochastic.
Throughout the paper assume that the routing matrix $R$ is out-connected, meaning that for every cell $i\in\mc E$ there exists some cell $j\in\mc E$ and an integer $k\ge0$ such that $\sum_{e\in\mc E}R_{je}<1$ and $(R^k)_{ij}>0$. 
Equation \eqref{eq:dyn1} can then be expressed more compactly as
\be\label{eq:dyn1bis}\dot{x} = \lambda - (I-R^T)z \,.\ee

In order to complete the description of the dynamical flow network system, it remains to specify how the outflow vector $z$ depends on the state vector $x$. In this paper, we focus on the case where the outflow $z_i$ from a cell $i$ is limited by a feedback-controller $\zeta_i(x)$, so that 
$$0\le z_i(t)\le \zeta_i(x(t))\,,\qquad i\in\mc E\,,\ t\ge0\,,$$
and that in fact the outflow $z_i$ from cell $i$ coincides with $\zeta_i(x)$ whenever the traffic volume $x_i$ is strictly positive, i.e., 
$$x_i(\zeta_i(x(t))-z_i)=0\,,\qquad i\in\mc E\,,\ t\ge0\,.$$
With this assumption, it is clear that the outflow from one link $z_i$ is only unspecified when there are no particles present on one link, i.e., when $x_i = 0$, while the controller still gives the link service such that $\zeta_i(x) > 0$. 
The rationale for not forcing $z_i=\zeta_i(x)$ also when $x_i=0$ is that in this case, if it happens also that $\lambda_i+\sum_jR_{ji}z_j<\zeta_i(x)$, then the dynamics would force the system to violate the non-negativity constraint on the traffic volume, i.e., 
$$x_i(t)\ge0\,,\qquad i\in\mc E\,,\ t\ge0\,.$$

Observe that the three constraints above may be rewritten more compactly as
\be\label{eq:dyn0} x\ge0\,,\ee
\be\label{eq:dyn2} 0 \leq z \leq \zeta(x)\,,\ee
\be\label{eq:dyn3}x^T(\zeta(x) - z) = 0 \, .\ee
where $\zeta:\mc X \to \R_+^\E$.  
Our main result presented as Theorem \ref{th:uniqueness} in Section \ref{sec:exsistence} guarantees that, whenver $\zeta$ is Lispchitz-continuous, the system of differential inclusions \eqref{eq:dyn1bis}--\eqref{eq:dyn3}  admits a unique solution for every initial state $x(0)\in\mc X$.

We conclude this section by discussing some examples in order to better motivate the considered dynamical network flow system and illustrate the usefulness of our result.

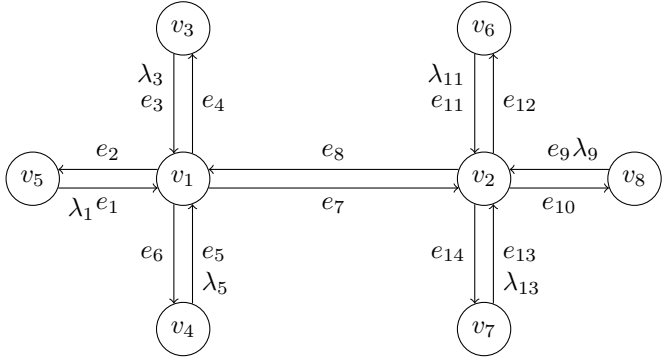
\begin{figure}
\centering
\begin{tikzpicture}
\node[draw,circle] (a) at (0,0) {$v_1$};
\node[circle,draw] (au) at (0,2) {$v_3$};
\node[circle,draw] (al) at (0,-2) {$v_4$};
\node[circle,draw] (aw) at (-2,0) {$v_5$};

\draw[->]  (aw.-20) node[below right] {$\lambda_1$} -- node[below] {$e_1$} (a.200);
\draw[->]  (a.160)  -- node[above] {$e_2$} (aw.20);

\draw[->] (au.-110) node[below left] {$\lambda_3$} -- node[left] {$e_3$} (a.110);
\draw[->] (a.70) -- node[right] {$e_4$} (au.-70);

\draw[->] (al.70) node[above right] {$\lambda_5$} -- node[right] {$e_5$} (a.-70);
\draw[->] (a.-110)  -- node[left] {$e_6$} (al.110);

\node[draw,circle] (b) at (4,0) {$v_2$};

\draw[->] (a.-20) -- node[below] {$e_7$} (b.200);
\draw[->] (b.160) -- node[above] {$e_8$} (a.20);

\node[circle, draw] (bu) at (4,2) {$v_6$};
\node[circle, draw] (bl) at (4,-2) {$v_7$};
\node[circle, draw] (be) at (6,0) {$v_8$};

\draw[->] (be.160) node[above left] {$\lambda_9$} -- node[above] {$e_9$} (b.20);
\draw[->] (b.-20)  -- node[below] {$e_{10}$} (be.200);

\draw[->] (bu.-110) node[below left] {$\lambda_{11}$} -- node[left] {$e_{11}$} (b.110);
\draw[->] (b.70)  -- node[right] {$e_{12}$} (bu.-70);

\draw[->] (bl.70) node[above right] {$\lambda_{13}$} -- node[right] {$e_{13}$} (b.-70);
\draw[->] (b.-110)  -- node[left] {$e_{14}$} (bl.110);

\end{tikzpicture}
\label{fig:twojunctions}
\caption {The network in Example~\ref{ex:signalizednetwork}. The network consists of two signalized junctions.}
\end{figure}

\begin{figure}
\centering
\input{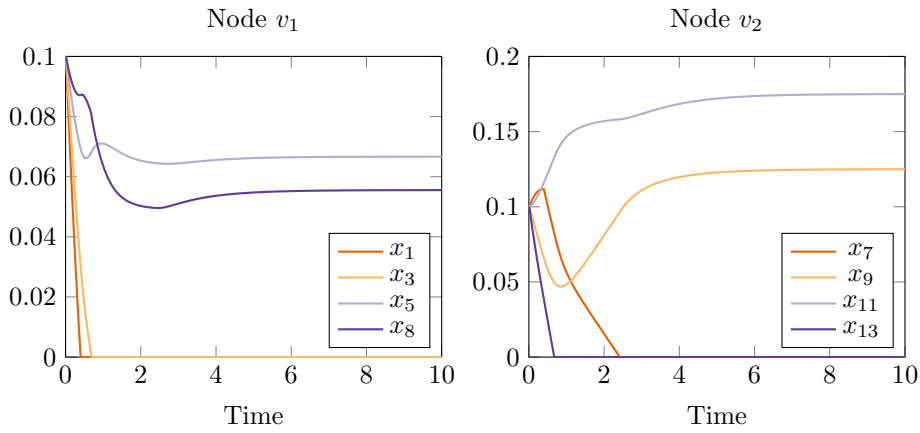}
\caption{How the traffic volumes $x$ evolve with time in Example~\ref{ex:signalizednetwork}.}
\label{fig:trafficnetworkvolumes}
\end{figure}

\begin{figure}
\centering
\input{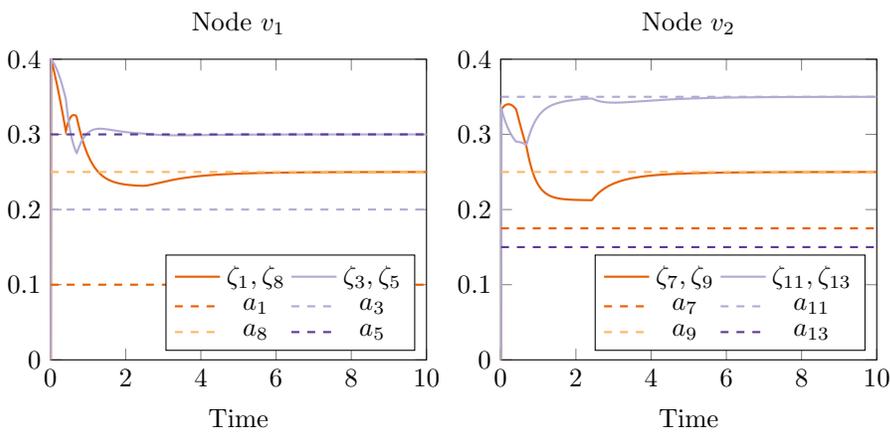}
\caption{How the control signals $\zeta$ evolve with time in Example~\ref{ex:signalizednetwork}. For reference we have also the equilibrium arrival rates $a = (I-R^T)^{-1}\lambda$ for each link. }
\label{fig:outflow}
\end{figure}

\begin{example}\label{ex:signalizednetwork}
Our framework applies to deterministic vertical queue network systems governed by feedback-controlled schedulers as used, e.g., in traffic signal control of urban networks \cite{Nilsson.Como:2019-TAC}.  To illustrate this application, consider a small dynamical flow network consisting of two controlled nodes as depicted in Figure~\ref{fig:twojunctions}. Suppose that nodes $v_1$ and $v_2$ are equipped with two service phases, such that either the east-west or north-south going links can receive service simultaneously. If the service allocation is computed by the Generalized Proportional Allocation (GPA) controller proposed in~\cite{nilsson2017generalized,Nilsson.Como:2019-TAC}, the outflow from each link will be limited by
$$\zeta_1 (x) = \zeta_8(x) =  \frac{x_1 + x_8}{x_1 + x_3 + x_5 + x_8 + \kappa_1} \,,$$ 
$$\zeta_3 (x) = \zeta_5(x) = \frac{x_3 + x_5}{x_1 + x_3 + x_5 + x_8 + \kappa_1}\,,$$
$$\zeta_7 (x) = \zeta_9(x) = \frac{x_7 + x_9}{x_7 + x_9 + x_{11} + x_{13} + \kappa_2} \,,$$ 
$$\zeta_{11} (x) = \zeta_{13}(x) = \frac{x_{11} + x_{13}}{x_7 + x_9 + x_{11} + x_{13} + \kappa_2}  \,,$$
where $\kappa_1, \kappa_2 > 0$ are constants, that are introduced to capture the fact that a fraction of the service cycle can not be utilized, due to the fact that there is some overhead time between the activation of subsequent phases.

For the links heading towards the boundary of the network, i.e., the links in the set $\mathcal B= \{e_{2}, e_{4}, e_{6}, e_{10},e_{12}, e_{14}\}$, we assume that particles are allowed to flow out from the network with unit capacity all times, i.e., $\zeta_i(x) =1$ for all $i \in \mathcal B$. Moreover, outflow from the boundary links will leave the network, so $R_{ij} = 0$ for all $i \in \mathcal B$ and all $j \in \mathcal E$.  

In this example, it is possible that control action is larger than the actual outflow. It can for instance happen when $x_1 =0$, but 
\begin{equation*}
\zeta_1(x) = \frac{x_8 }{ x_3 + x_5 + x_8  + \kappa_1} > \lambda_1 \,.
 \end{equation*}
A simulation of the dynamical flow network is shown in Figure~\ref{fig:trafficnetworkvolumes}. In the simulation, suppose that $1/4$ of the inflow from each link to the nodes $v_1$ and $v_2$ is going left, $1/4$ going right, and the remaining half of the flow proceeds straight. Moreover, we let $\lambda_1 = 0.10$, $\lambda_3 = 0.20$, $\lambda_5 = 0.30$, $\lambda_9 = 0.25$, $\lambda_{11} = 0.35$, and $\lambda_{13} = 0.15$. The constants in the controllers are chosen to be $\kappa_1 = 0.1$ and $\kappa_2 = 0.2$, and all the traffic volumes are initiated at $0.1$, i.e., $x_i(0) = 0.1$ for all $i \in \mathcal{E}$. 

In Figure~\ref{fig:outflow} the control actions $\zeta_i(x)$ are plotted, together with the outflows at equilibrium $a$. The latter can be computed through $a = (I-R^T)^{-1}\lambda$. From Figure~\ref{fig:trafficnetworkvolumes} and Figure~\ref{fig:outflow}, we can see that controller will be equal to the equilibrium flows for all links where $x_i > 0$, while for the links where $x_i = 0$, the controller will allow for more outflow than what is physically possible, and hence $z_i < \zeta_i(x_i)$ for those links. Moreover, the control action converges to the outflows for the links with $x_i > 0$. This observation, is a consequence of the fact proven in~\cite{nilsson2017generalized} that the GPA controller is stabilizing.
\end{example}

While the example above assume no propagation delay between the nodes, dynamical propagation delay can be introduced in model by adding intermediate nodes, as the following example shows. This property makes the proposed model more adaptable to certain applications, compared to the open loop model presented in~\cite{muralidharan2015analysis, hosseini2017}, where the delay is assumed to be independent of the state.

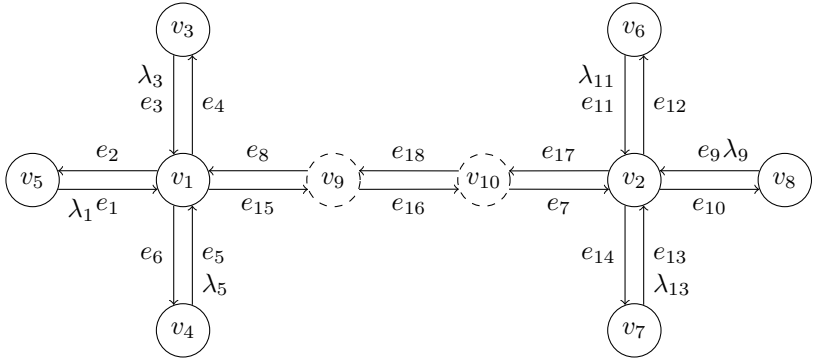
\begin{figure}
\centering
\begin{tikzpicture}
\node[draw,circle ] (a) at (0,0) {$v_1$};
\node[circle,draw ] (au) at (0,2) {$v_3$};
\node[circle, draw] (al) at (0,-2) {$v_4$};
\node[circle, draw] (aw) at (-2,0) {$v_5$};

\draw[->]  (aw.-20) node[below right] {$\lambda_1$} -- node[below] {$e_1$} (a.200);
\draw[->]  (a.160)  -- node[above] {$e_2$} (aw.20);

\draw[->] (au.-110) node[below left] {$\lambda_3$} -- node[left] {$e_3$} (a.110);
\draw[->] (a.70) -- node[right] {$e_4$} (au.-70);

\draw[->] (al.70) node[above right] {$\lambda_5$} -- node[right] {$e_5$} (a.-70);
\draw[->] (a.-110)  -- node[left] {$e_6$} (al.110);

\node[draw,circle] (b) at (6,0) {$v_2$};
\node[draw,dashed,circle] (v3) at (2,0) {$v_9$};
\node[draw,dashed,circle, inner sep=2pt] (v4) at (4,0) {$v_{10}$};

\draw[->] (v4.-20) -- node[below] {$e_7$} (b.200);
\draw[->] (v3.160) -- node[above] {$e_8$} (a.20);

\draw[->] (a.-20) -- node[below] {$e_{15}$} (v3.200);
\draw[->] (v3.-20) -- node[below] {$e_{16}$} (v4.200);

\draw[->] (b.160) -- node[above] {$e_{17}$} (v4.20);
\draw[->] (v4.160) -- node[above] {$e_{18}$} (v3.20);

\node[circle, draw] (bu) at (6,2) {$v_6$};
\node[circle, draw] (bl) at (6,-2) {$v_7$};
\node[circle, draw] (be) at (8,0) {$v_8$};

\draw[->] (be.160) node[above left] {$\lambda_9$} -- node[above] {$e_9$} (b.20);
\draw[->] (b.-20)  -- node[below] {$e_{10}$} (be.200);

\draw[->] (bu.-110) node[below left] {$\lambda_{11}$} -- node[left] {$e_{11}$} (b.110);
\draw[->] (b.70)  -- node[right] {$e_{12}$} (bu.-70);

\draw[->] (bl.70) node[above right] {$\lambda_{13}$} -- node[right] {$e_{13}$} (b.-70);
\draw[->] (b.-110)  -- node[left] {$e_{14}$} (bl.110);

\end{tikzpicture}
\caption {The network in Example~\ref{ex:ctm}. By introducing intermediate nodes between the junction, the flow dynamics can be discretized and a dynamic propagation delay can be modeled. }
\label{fig:networkctm}
\end{figure}

\begin{figure}
\centering
\input{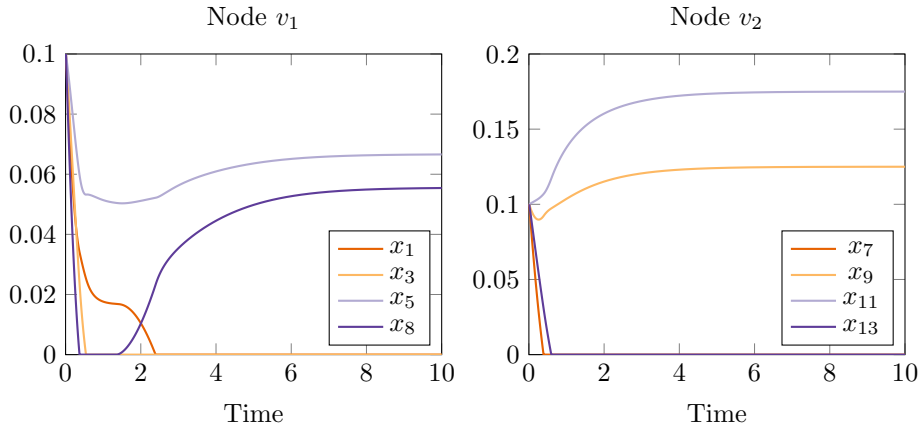}
\caption{How the traffic volumes $x$ evolve with time in Example~\ref{ex:ctm}. Compared to Example~\ref{ex:signalizednetwork}, the trajectories becomes different, due to the flow propagation dynamics on the intermediate links.}
\label{fig:ctmtra}
\end{figure}

\begin{example}\label{ex:ctm}
Starting from Example~\ref{ex:signalizednetwork}, we introduce two intermediate nodes between the nodes $v_1$, and $v_2$, as shown in Figure ~~\ref{fig:networkctm}. Moreover, we let the outflow from the added intermediate links adhere the continuous version of the Cell Transmission Model  (CTM)~\cite{Daganzo:94, Daganzo:95}, used to model traffic flow in e.g.~\cite{lovisari2014monotone,coogan2018polytree}. To each of the links $e_{15}, e_{16}, e_{17}$ and $e_{18}$ we assign a demand function $d_i(x_i) \geq 0$ that is strictly increasing and Lipschitz continuous in the traffic volume. To each of links $e_7, e_8, e_{16}$ and $e_{18}$ we assign a supply function $s_i(x_i) \geq 0$ that is non-increasing and Lipschitz continuous in the traffic volume. The outflows from the the intermediate links are then given by
\begin{align*}
\zeta_{15}(x) &= \min(d_{15}(x_{15}), s_{16}(x_{16}))\,, & \zeta_{16}(x) &= \min(d_{16}(x_{16}), s_7(x_7)) \,, \\
\zeta_{17}(x) &= \min(d_{17}(x_{17}), s_{18}(x_{18}))\,, & \zeta_{18}(x) &= \min(d_{18}(x_{18}), s_7(x_8)) \,. \\
\end{align*}

In Figure~\ref{fig:ctmtra} we show the trajectories for the traffic volumes on the incoming links to node $v_1$ and $v_2$. The simulation parameters are the same as in Example~\ref{ex:signalizednetwork}, and for all the intermediate links we let $d_i(x_i) = x_i$ and
\begin{equation*}
s_i(x_i) = \begin{cases} 
1-x_i & \text{if } x_i \leq 1 \,, \\
0 & \text{otherwise.}
\end{cases}
\end{equation*}
Moreover, we let $x_{15} (0) = x_{16} (0) = x_{17}(0) = x_{18} (0) = 0$. That the intermediate nodes introduces a propagation delay, can be seen though that it takes a longer time for the traffic volume on link $e_8$ to converge. The reason for that the same delay can not be observed on link $e_7$, is that the controller is already allowing for more outflow than needed, due to a high traffic volume on link $e_9$.

As a remark, in the case when the demand functions are on the form
\begin{align*}
d_i(x_i) = C_i \frac{x_i}{x_i + \kappa_i} \,,
\end{align*}
where $C_i >0$ and $\kappa_i > 0$ are constants. If $s_j(x_j) \geq C_i$ for all $x_j$, then the stability analysis for the General Proportional Allocation controller, done in~\cite{nilsson2017generalized}, can be applied to ensure stability of the dynamical flow network. This since the demand function is in fact a GPA controller with just one incoming link.
\end{example}

\section{Existence and Uniqueness}\label{sec:exsistence}
In this section, we present a proof of existence and uniqueness of a solution to the dynamical system~\eqref{eq:dyn1bis}--\eqref{eq:dyn3}. The proof is an extension of the reflection principle for Brownian motion, previously presented in~\cite{harrison1981reflected}, to the case when the outflow is determined by a Lipschitz continuous controller.

\begin{theorem}\label{th:uniqueness}
Let $R$ be an out-connected routing matrix and $\lambda:\R_+\to \R_+^\mc E$ is a bounded measurable function depending on time. Assume that $\zeta:\mc X\to\mc X$ Lipschitz continuous, then for every initial condition $x(0) \in \mc X$, the dynamics given by~\eqref{eq:dyn1bis}--\eqref{eq:dyn3} admits a unique solution.
\end{theorem}
Throughout the proof the of Theorem \ref{th:uniqueness}, we will make use of the fact that since the routing matrix is out-connected, it has a spectral radius strictly smaller than $1$, see e.g.~\cite[Proof of Theorem 2]{como2016local}.  Then, from~\cite[Lemma 5.6.10]{horn2012matrix} it follows that there exists a norm $\norm{\cdot}_\dagger$ on $\R^{\mc E}$ such that the induced matrix norm of $R$ satisfies $\norm{R}_\dagger < 1$. For $T>0$, we shall consider the space $\mc C_{T}$ of continuous vector-valued functions $f:[0,T]\to\R^{\mc E}$ equipped with the norm  
$$\norm{f} = \norm{\sup_{0 \leq t \leq T} |f(t)| }_\dagger\, .$$

As a first step towards the proof of Theorem \ref{th:uniqueness}, to a given continuous vector-valued function $\gamma$ in $\mc C_T$, we associate the operator $\Pi_{\gamma}:\mc C_T\to\mc C_T$ defined, for  $0\le t\le T$, by 
\be\label{Pi-def}{\left[\Pi_{\gamma}(v)\right](t)} = \sup_{0 \leq s \leq t} \left[ R^Tv(s) - \gamma(s)\right]_+\,.\ee The operator $\Pi_{\gamma}$ has the following fundamental properties.

\begin{lemma}\label{lemma:fixedpoint} For every $T>0$ and every continuous vector-valued function $\gamma$ in $\mc C_T$, the operator
$\Pi_{\gamma}$ is a contraction  on $\mc C_T$.
\end{lemma}
\proof
We will first show that $\Pi_{\gamma}$ is a contraction on $\mc C_T$. For any $v, w$ in $\mc C_T$,  $0\le s\le  t\le T$, and 
$i$ in $\mc E$, put \be\label{fg-def}f(s)=[R^T v(s) - \gamma(s)]_i\,,\quad g(s)=[R^T w(s) - \gamma(s)]_i\,,\ee
$$h(s)=f(s)-g(s)\,.$$
Choose some 
$$s^* \in\arg\max_{0 \leq s \leq t}[f(s)]_+\,,\qquad q^* \in \arg\max_{0 \leq s \leq t}[g(s)]_+\,,$$
and observe that 
\be\label{+1}[f(s^*)]_+=[g(s^*)+h(s^*))]_+\le[g(s^*)]_++[h(s^*))]_+\,,\ee
\be\label{+2}[g(q^*)]_+=[f(q^*)-h(q^*))]_+\le[f(q^*)]_++[-h(q^*))]_+=[f(q^*)]_++[h(q^*))]_-\,.\ee
Using \eqref{+1} and the fact that $[f(s^*)]_+=\ds\sup_{0 \leq s \leq t} [f(s)]_+$, we get
$$
\ba{rcl}
\ds\sup_{0 \leq s \leq t}[h(s)]_+
&\ge&\ds[h(s^*)]_+\\[10pt]
&\ge&\ds[f(s^*)]_+ - [g(s^*)]_+\\[10pt] 
&\ge&\ds\sup_{0 \leq s \leq t} [f(s)]_+ - \sup_{0 \leq s \leq t} [g(s)]_+\,.\ea
$$
Analogously, \eqref{+2} and the fact that $[g(q^*)]_+=\ds\sup_{0 \leq s \leq t} [g(s)]_+$ give
$$\ba{rcl}
\ds\sup_{0 \leq s \leq t}[h(s)]_-  
&=&\ds  \sup_{0 \leq s \leq t}[-h(s)]_+ \\[10pt]
& \ge&\ds [g(q^*)]_+ - [f(q^*)]_+ \\[10pt]
&\ge&\ds  \sup_{0 \leq s \leq t} [g(s)]_+- \sup_{0 \leq s \leq t} [f(s)]_+ \,.
\ea$$
Therefore, 
$$\ba{rcl}\ds\sup_{0 \leq s \leq t} |h(s)|
&=&
\ds\max\left\{\sup_{0 \leq s \leq t} [h(s)]_+,\sup_{0 \leq s \leq t}  [h(s)]_-\right\}\\
&\geq&\ds\left|\sup_{0 \leq s \leq t} [f(s)]_+-\sup_{0 \leq s \leq t}  [g(s)]_+\right|\,. 
\ea$$  
Now, let us define the vector $\alpha\in\R^{\mc E}$ with entries 
$$ \alpha_i=\sup_{0 \leq t \leq T}\left|\left[\Pi_{\gamma}(v)\right]_i(t)-\left[\Pi_{\gamma}(w)\right]_i(t)\right|\,,$$
for all $i\in\mc E$. Using \eqref{Pi-def} and \eqref{fg-def}, we get
$$\ba{rcl}\ds\alpha_i
&=&\ds\sup_{0 \leq t \leq T}\left|\sup_{0 \leq s \leq t}\left[ f(s)\right]_+-\sup_{0 \leq s \leq t} \left[ g(s)\right]_+\right|\\[10pt] 
&\le&\ds\sup_{0 \leq t \leq T}\sup_{0 \leq s \leq t} \left| h(s)\right|\\[10pt] 
&=&\ds\sup_{0 \leq t \leq T} \left|[R^T(v(t)-w(t))]_i\right|\\[10pt] 
&\le&\ds \sum_jR_{ji}\sup_{0 \leq t \leq T}|v_j(t)-w_j(t)|\,.\ea
$$
Hence, 
$$
||\Pi_{\gamma}v-\Pi_{\gamma}w||=||\alpha||_\dagger\le||R^T||_\dagger||v-w||
$$
Since $||R^T||_\dagger< 1$, the above proves that $\Pi_{\gamma}$ is a contraction on $\mc C_T$.\qed

It follows from Lemma \ref{lemma:fixedpoint} and the Banach fixed point theorem that  for every continuous vector-valued function $\gamma$ in $\mc C_T$, the operator
$\Pi_{\gamma}$ admits a unique fixed point 
\be\label{Psi-def}\Psi(\gamma)= \Pi_{\gamma}(\Psi(\gamma))\in\mc C_T\,.\ee 
The following result characterizes the dependance of such fixed point $\Psi(\gamma)$ on the choice of $\gamma$ in $\mc C_T$. 

\begin{lemma}\label{lemma:fixedpoint-bis}
For every $T>0$, the operator $\Psi:\mc C_T\to\mc C_T$
that maps a continuous vector-valued function $\gamma$ into the unique fixed point of the associated operator $\Pi_{\gamma}$ is Lipschitz continuous.
\end{lemma}
\proof
For $k\ge0$, let $\Pi_{\gamma}^{k}$ be the composition of $\Pi_{\gamma}$ with itself $k$ times. Fix three functions $v,\gamma, \eta \in \mc C_T$ and for $0\le t\le T$, define 
$$\delta^k(t)=\left[\Pi_{\gamma}^{k}(v)\right](t) - \left[\Pi_{\eta}^{k}(v)\right](t)\,.$$
Then, we have that 
$$\ba{rcl}
\left|\delta^{k+1}(t)\right|
&=&\left|\left[\Pi_{\gamma}^{k+1}(v)\right](t) - \left[\Pi_{\eta}^{k+1}(v)\right](t)\right|\\[10pt]
&=&\ds\left|\sup_{0 \leq s \leq t} \left[ R^T[\Pi_{\gamma}^{k}v](s) - \gamma(s) \right]_+ \right. \left. - \sup_{0 \leq s \leq t} \left[ R^T [\Pi_{\eta}^{k}v](s) - \eta(s) \right]_+\right|  \\[10pt]
&\leq&\ds\left|\sup_{0 \leq s \leq t} \left[ R^T\left([\Pi_{\gamma}^{k}v](s) - [\Pi_{\eta}^{k}v](s)\right) - \left(\gamma(s) - \eta(s) \right) \right]_+ \right| \\[10pt]
&\leq& \ds \sup_{0 \leq s \leq t} \left| R^T\delta^k(s) \right| +\sup_{0 \leq s \leq t}|\gamma(s) - \eta(s)|\\[10pt] 
\ea$$
so that 
$$
||\delta^{k+1}||\le||R^T||_\dagger||\delta^{k}||+||\gamma-\eta||\,.
$$
It follows that, for all $v$ in $\mc C_T$ and $k\ge0$, 
$$\norm{\Pi_{\gamma}^{k}(v)- \Pi_{\eta}^{k}(v)}\leq \sum_{l = 0}^k ||R^T||_\dagger^l \norm{\gamma- \eta} \,.$$ 
Since $||R^T||_\dagger<1$ and $\Pi_{\gamma}$ and $\Pi_{\eta}$ are both contractions with fixed points $\Psi(\gamma)$ and $\Psi(\eta)$, respectively, taking the limit as $k$ grows large in the above gives
$$\ba{rclcl}\ds\norm{\Psi(\gamma)- \Psi(\eta)}
&=&\ds\lim_{k\to\infty}\norm{\Pi_{\gamma}^{k}(v)- \Pi_{\eta}^{k}(v)}\\[10pt] 
&\leq&\ds \sum_{l = 0}^{+\infty} ||R^T||_\dagger^l \norm{\gamma- \eta}\\[10pt]
&=&\ds\frac{\norm{\gamma- \eta}}{1-||R^T||_\dagger} \,.\ea$$ 
which concludes the proof of the Lemma~\ref{lemma:fixedpoint}.
\qed

Our next step towards proving Theorem \ref{th:uniqueness} consists in finding an equivalent formulation of the controlled traffic network dynamics \eqref{eq:dyn1bis}--\eqref{eq:dyn3}.
Towards this goal, we introduce two operators $$\Phi,\Gamma:\mc C_T\to\mc C_T$$ 
defined by 
\be\label{Phi-def}\Phi(y)=y + (I -R^T)\Psi(y)\ee
and, respectively,
\be\label{Gamma-def}\Gamma(x)(t) =  x(0) + \int_0^t \left( \lambda(s) - (I -R^T) \zeta(x(s))  \right) \de s \, .\ee
The relationship between the controlled traffic network dynamics \eqref{eq:dyn1bis}--\eqref{eq:dyn3} and the operators \eqref{Phi-def} and \eqref{Gamma-def} is illustrated by the following results. 

\begin{lemma}\label{lemma:new-equivalence}
Let $R$ be an out-connected routing matrix, $\lambda$ a possible time varying exogenous inflow vector, $\zeta:\mc X\to\mc Z$ a Lipschitz continuous function, and  $x(0) \in \mc X$. Then, $(x(t),z(t))$ is a solution of the controlled traffic network dynamics \eqref{eq:dyn1bis}--\eqref{eq:dyn3} in a time interval $[0,T]$ with initial condition $x(0)$ if and only if there exist $y,w\in\mc C_T$ that are absolutely continuous and such that 
\be\label{x=Phi(y)}x=\Phi(y)\,,\ee
\be\label{y=Gamma(x)}y=\Gamma(x)\,,\ee
\be\label{w=Psi(y)}w=\Psi(y)\,,\ee
and, for almost all $t\in[0,T]$,
\be\label{z=zeta(x)-wdot}z=\zeta(x)-\dot w\,.\ee
\end{lemma}

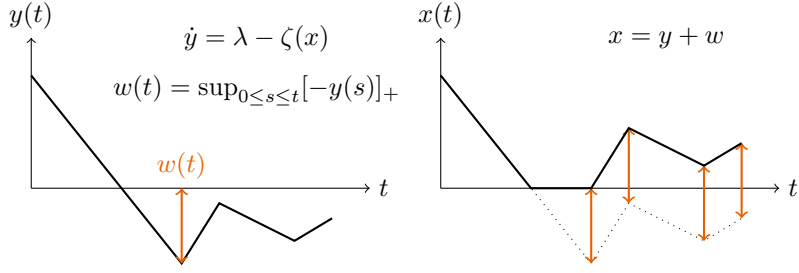
\begin{figure}
\centering
\begin{tikzpicture}
\draw[->] (0,0) -- (0,2) node[above] {$y(t)$};
\draw[->] (0,0) -- (4.5,0) node[right] {$t$};

\draw[thick] (0,1.5) -- (2, -1) -- (2.5,-0.2) -- (3.5, -0.7) -- (4, -0.4) ;

\draw[mycolor1, thick, <->]  (2, -1) -- (2,0) node[above]{$w(t)$} ;

\node (0) at (3, 2) {$\dot{y} = \lambda - \zeta(x)$};
\node (0) at (3, 1.3) {$w(t) = \sup_{0 \leq s \leq t} [-y(s)]_+$};
\end{tikzpicture}
\begin{tikzpicture}
\draw[->] (0,0) -- (0,2) node[above] {$x(t)$};
\draw[->] (0,0) -- (4.5,0) node[right] {$t$};

\draw[dotted] (0,1.5) -- (2, -1) -- (2.5,-0.2) -- (3.5, -0.7) -- (4, -0.4) ;

\draw[mycolor1, thick, <->]  (2, -1) -- +(0,1);
\draw[mycolor1, thick, <->]  (2.5,-0.2) -- +(0,1);
\draw[mycolor1, thick, <->]  (3.5, -0.7)  -- +(0,1);
\draw[mycolor1, thick, <->]  (4, -0.4) -- +(0,1);

\draw[thick] (0,1.5) -- (1.2,0) -- (2, 0) -- (2.5,0.8) -- (3.5, 0.3) -- (4, 0.6)  ;

\node (0) at (3, 2) {$x = y + w$};
\end{tikzpicture}
\caption{The connection between the quantities $x$, $y$ and $w$ in Lemma~\ref{lemma:new-equivalence} for the case when the network consists of a single cell}
\label{fig:ywxexplained}
\end{figure}

\begin{remark}
There is an interpretation of the quantities in Lemma~\ref{lemma:new-equivalence}. The quantity $y$ can be seen as the traffic volumes on the links \emph{if} the volumes were allowed to be negative, and $w$ is how much one must add to this quantity to make sure that the traffic volume $x$ stays non-negative. In Figure~\ref{fig:ywxexplained} those trajectories are illustrated for a single cell, i.e., $R = 0$. Observe that $w(t)$ is non-decreasing and only increases when $x = 0$.
\end{remark}

\proof
(i) Let $(x(t),z(t))$ be a solution of the controlled traffic network dynamics \eqref{eq:dyn1bis}--\eqref{eq:dyn3} on $[0,T]$ with initial condition $x(0)$. 
For $0\le t\le T$, let 
\be\label{w}w(t)=\int_0^t(\zeta(x(s)) - z(s))\de s\,,\ee
\be\label{y}y(t)=x(t)-(I - R^T)w(t)\,.\ee
We will show that \eqref{x=Phi(y)}--\eqref{z=zeta(x)-wdot} are satisfied. 
Indeed, taking the time derivative of both sides of \eqref{w} gives \eqref{z=zeta(x)-wdot}. 
On the other hand, \eqref{y}, \eqref{eq:dyn1bis}, and \eqref{w} yield
$$\ba{rcl} y(t)&=& x(t) - (I-R^T) w(t) \\[10pt] 
&=&\ds x(0)+\int_0^t(\lambda(s)-(I-R^T)z(s))\de s   -(I - R^T)w(t)\\[10pt] 
&=&\ds x(0)+\int_0^t(\lambda(s)-(I-R^T)\zeta(x(s)))\de s\\[10pt]
&=&\Gamma(x)(t)\,,
\ea $$
so that \eqref{y=Gamma(x)} is satisfied as well. Moreover, \eqref{y} and \eqref{w=Psi(y)} clearly imply \eqref{x=Phi(y)}. Hence, it remains to prove \eqref{w=Psi(y)}.
For that, first observe that \eqref{w} and \eqref{eq:dyn3} imply that 
\be \label{eq:equivdyn3} 
x\ge0\,,\qquad x^T \dot{w} = 0 \,, \qquad 0 \leq \dot{w} \leq \zeta(x)  \,. \ee
In turn, the above and \eqref{y} imply that,  for $0\le s\le t$, 
$$w(t)\ge w(s) = R^Tw(s) + x(s) - y(s)\ge R^Tw(s) - y(s)\, ,$$
so that 
$$w(t) \geq \sup_{0 \leq s \leq t}\left\{R^T w(s) - y(s)\right\}\, .$$
Since $w(0)$ is non-increasing and $w(0)=0$, we have $w(t)\ge0$, which together with the above gives 
$$w(t) \geq \sup_{0 \leq s \leq t} \left[ R^T w(s) - y(s) \right]_+=\Pi_{y}(w)(t)  \, .$$
In fact, if the above were not an identity for some $0\le t\le T$, there would exist some $0\le t^*\le T$ and $i\in\mc E$ such that 
\be\label{contradiction}w_i(t^*)>\sup_{0 \leq s \leq t^*} \left\{\sum\nolimits_jR_{ji}w_j(s) - y_i(s)\right\}\,,\qquad \dot w_i(t^*)>0\,.\ee
But the second inequality above and  \eqref{eq:equivdyn3} imply that $x_i(t^*)=0$ so that, by \eqref{y}, 
$y_i(t^*)=\sum\nolimits_jR_{ji}w_j(t^*) - y_i(t^*)$ which contradicts \eqref{contradiction}. Hence, we necessarily have 
$$w(t) =\Pi_{y}(w)(t)\,,\qquad 0\le t\le T \,,$$
i.e., $w$ is the fixed point $\Pi_y$ on $\mc C_T$, so that \eqref{w=Psi(y)} is satisfied.

(ii) Let $w,x,y,z\in\mc C_T$ be such that $y$ and $w$ are absolutely continuous and \eqref{x=Phi(y)}--\eqref{z=zeta(x)-wdot} are satisfied. Then, for $0\le t\le T$, an application of \eqref{x=Phi(y)}, \eqref{Phi-def}, \eqref{y=Gamma(x)}, \eqref{w=Psi(y)}, \eqref{Gamma-def}, and \eqref{z=zeta(x)-wdot} give 
$$\ba{rcl}x(t)
&=&\Phi(y)(t)\\[10pt]
&=&y(t)+(I-R^T)\Psi(y)(t)\\[10pt]
&=&\Gamma(x)(t)+(I-R^T)w(t)\\[10pt]
&=&\ds x(0) + \int_0^t \!\left( \lambda(s) - (I -R^T) \zeta(x(s))  \right) \de s +  (I -R^T)\int_0^t \left( \zeta(x(s))-z(s) \right) \de s\\[10pt] 
&=&\ds x(0) + \int_0^t \left( \lambda(s) - (I -R^T) z(s)  \right) \de s\,,
\ea
$$
hence \eqref{eq:dyn1bis} is satisfied. 
On the other hand, \eqref{x=Phi(y)}, \eqref{Phi-def}, \eqref{Psi-def}, and \eqref{Pi-def} give 
\be\label{x-constr}\ba{rcl}x(t)
&=&\ds \Phi(y)(t)\\[10pt]
&=&\ds y(t)+(I-R^T)\Psi(y)(t)\\[10pt]
&=&\ds y(t)-R^T\Psi(y)(t)+\sup_{0\le s\le t}\left[R^T\Psi(y)(s)-y(s)\right]_+\\[10pt]
&\ge&\ds0\,. 
\ea\ee
Moreover, 
\eqref{w=Psi(y)}, \eqref{Psi-def}, and \eqref{Pi-def} yield 
\be\label{w-bis}w(t)=\Psi(y)(t)=\sup_{0 \leq s \leq t} \left[ R^Tw(s) - y(s)\right]_+\,,\ee
so that $w_i(t)$ is non-decreasing for all $i\in\mc E$, hence $\dot{w} \geq 0\,.$ Furthermore, let $\mc I:=\{i\in\mc E:\,\dot w_i(t)>0\}$ be the set of cells $i$ such that $w_i(t)$ is strictly increasing at time $t$. It then follows from  \eqref{w-bis} that 
\be\label{wi}w_i(t) = \sum_{j\in\mc E}R_{ji}w_j(t) -y_i(t)\,,\qquad i\in\mc I\,.\ee
Equation \eqref{wi} implies that, for $i\in\mc I$, 
$$\ba{rcl}
\ds\dot w_i(t)
&=&\ds\sum_{j\in\mc E}R_{ji}\dot w_j(t) -\dot y_i(t)\\[10pt]
&=&\ds\sum_{j\in\mc I}R_{ji}\dot w_j(t) -\lambda_i(t)+\zeta_i(x(t))-\sum_{j\in\mc E}R_{ji}\zeta_j(x(t))\\[10pt] 
&\le&\ds\sum_{j\in\mc I}R_{ji}\dot w_j(t) -\lambda_i(t)+\zeta_i(x(t))-\sum_{j\in\mc I}R_{ji}\zeta_j(x(t))
\,.\ea$$
The above implies that  
\be\label{I-RII}(I-R_{\mc I\mc I}^T)\dot w_{\mc I}(t)\le (I-R_{\mc I\mc I}^T)\zeta_{\mc I}(x(t))-\lambda_{\mc I}(t)\,,\ee
where $R_{\mc I\mc I}$ is the $\mc I\times\mc I$ block of $R$ and $\dot w_{\mc I}(t)$, $\zeta_{\mc I}(x(t))$, and $\lambda_{\mc I}(t)$ are the $\mc I$ blocks of the corresponding vectors $\dot w(t)$, $\zeta_{\mc I}(x(t))$, and $\lambda_{\mc I}(t)$. Since $R$ is out-connected, each of its diagonal blocks such as $R_{\mc I\mc I}$   has spectral radius smaller than $1$. Hence $(I-R_{\mc I\mc I}^T)$ invertible with nonnegative inverse $(I-R_{\mc I\mc I}^T)^{-1}$. Hence, \eqref{I-RII} implies that 
$$\dot w_{\mc I}(t)\le\zeta_{\mc I}(x(t))-(I-R_{\mc I\mc I}^T)^{-1}\lambda_{\mc I}(t)\le\zeta_{\mc I}(x(t))\,.$$
Since $\dot w_{\mc E\setminus\mc I}(t)=0$ by definition and we have already noticed that $\dot w(t)\ge0$, we thus have that $z=\zeta(x)-\dot w$ satisfies
\be\label{z-constr}0\le z\le\zeta(x)\,.\ee 
Finally, using again \eqref{x=Phi(y)}, \eqref{Phi-def}, \eqref{w=Psi(y)}, and \eqref{wi}, one gets that 
$$x_i(t)=y_i(t)+w_i(t)-\sum_jR_{ji}w_j(t)=0$$
for every $i\in\mc I$. Along with \eqref{x-constr} and \eqref{z-constr}, this implies that 
\be\label{xz-constr}x^T(\zeta(x)-z)=x^T\dot w=0\,.\ee
From \eqref{x-constr}, \eqref{z-constr}, and \eqref{xz-constr} it follows that \eqref{eq:dyn3} is satisfied. Therefore $(x,z)$ is a solution of \eqref{eq:dyn1bis}--\eqref{eq:dyn3}. 
\qed

We are now in a position to prove Theorem~\ref{th:uniqueness}. It follows from Lemma~\ref{lemma:fixedpoint-bis} that $\Psi$ is a Lipschitz continuous operator on $\mc C_T$.  Hence, the operator $\Phi$ is Lipschitz-continuous as well and we shall denote by $\phi > 0$ its Lipschitz constant.
Since $\zeta:\mc X\to\R^{\mc E}$ is a Lipschitz continuous function, the operator $\Gamma$ is Lipschitz-continuous on $\mc C_T$ for all $T > 0$, with Lipschitz constant equal to $\varpi T$ for some constant $\varpi > 0$ that is independent from $T$. It then follows that, for  $0 < T < (\varpi \phi)^{-1}$, the composition operator $\Phi \circ \Gamma:\mc C_T\to\mc C_T$ is Lipschitz continuous with Lipchitz constant $$L = \varpi \phi T < 1\,.$$ Therefore, $\Phi \circ \Gamma$ is a contraction on $\mc  C_T $, hence it has a unique fixed point $x = \Phi(\Gamma(x))$. Let $$y = \Gamma(x)\,,\qquad w = \Psi(y)\,,\qquad z=\zeta(x)-\dot w\,.$$ 
By Lemma~\ref{lemma:new-equivalence} we get that this $(x,z)$ is the unique solution to~\eqref{eq:dyn1bis}--\eqref{eq:dyn3} on $[0,T]$ with initial condition $x(0)$. Existence and uniqueness of the solution $(x,z)$ of~\eqref{eq:dyn1bis}--\eqref{eq:dyn3} can then be extended to the whole semi-infinite time interval $[0,+\infty)$ by standard continuation arguments.

\section{Conclusions}
In this paper, we have presented a dynamical flow network model, where a feedback-controller limits the outflow from each link. The feedback-controller may allow for more outflow that is physically possible to flow out. Due to this property, the dynamical system is described as a differential inclusion. We show that such differential inclusion admits a unique solution. In the future, we plan to extend the well-posedness results, at least for the existence part, to non-Lipschitz and possibly discontinuous feedback controls as those mentioned in \cite{Nilsson.Como:2019-TAC}. It would also be of great interest to study the case of time-varying routing matrix $R$ and/or to analyze how robust the model is to the choice of such matrix.

\bibliographystyle{siamplain}
\bibliography{ref}             

\end{document}